\documentclass[12pt,a4paper]{article}%
\usepackage[utf8]{inputenc}
\usepackage{hyperref}
\usepackage{amsmath}
\usepackage{amssymb}
\newtheorem{theorem}{Theorem}
\newtheorem{conjecture}[theorem]{Conjecture}
\newtheorem{lemma}[theorem]{Lemma}
\newtheorem{proposition}[theorem]{Proposition}
\newenvironment{proof}[1][Proof]{\noindent\textbf{#1.} }{\ \hfill \rule{0.5em}{0.5em}\bigskip}
\begin{document}

\hypersetup{pdftitle={The Manhattan cospectrality of the cocktail party graph}}
\title{The Manhattan cospectrality of the cocktail party graph}
\author{Irena M. Jovanovi\'{c}$^{1}$, Jelena Sedlar$^{2,4}$, Riste
\v{S}krekovski$^{3,4}$\\[0.3cm]
{\small $^{1}$ \textit{Independent Researcher, Belgrade, Serbia}}\\[0.1cm]
{\small $^{2}$ \textit{University of Split, Faculty of Civil Engineering,
Architecture and Geodesy, Croatia}}\\[0.1cm]
{\small $^{3}$ \textit{Faculty of Mathematics and Physics, University of
Ljubljana, Slovenia}}\\[0.1cm]
{\small $^{4}$ \textit{Faculty of Information Studies, Novo mesto, Slovenia}}}
\date{}
\maketitle

\begin{abstract}
For a graph $G$ of order $N$ with adjacency eigenvalues $\lambda_{1}(G)\geq
\cdots\geq\lambda_{N}(G),$ the cospectrality of $G$ is the smallest value of
$\sum_{i}\left\vert \lambda_{i}(G)-\lambda_{i}(H)\right\vert $ over all graphs
$H$ of order $N$ not isomorphic to $G.$ Abdollahi and Zakeri (Ars Math. Contemp. 22(1) (2022) P1.06) conjectured that
the cospectrality of the cocktail party graph $CP_{n}$ is attained exactly at
the graph $H_{n}.$ This graph is obtained from $K_{1,2,\ldots,2}$ by attaching
a pendant vertex to the vertex of the part of size one. We prove this
conjecture for every $n\geq2$ by showing that only three explicit graphs (two
for $n=2$) are at distance less than $2$ from $CP_{n}.$ We work with the
complement of a competing graph. First we remove its components $K_{2}.$
Then we treat the case in which its main component has a dominating vertex.
The remaining cases are reduced to finitely many small graphs. The proof is
computer-assisted, and all computations concern explicit finite lists of
graphs.
\end{abstract}

\section{Introduction}

All graphs here are finite, simple and undirected. For a graph $G$ of order
$N$ we write $\lambda_{1}(G)\geq\lambda_{2}(G)\geq\cdots\geq\lambda_{N}(G)$
for the eigenvalues of its adjacency matrix $A(G).$ Let $G$ and $H$ be graphs
of the same order $N.$ We compare their eigenvalues in non-increasing order,
position by position. The (Manhattan) \emph{spectral distance} of $G$ and $H$
and the \emph{cospectrality} of $G$ are
\[
\sigma(G,H)=\sum_{i=1}^{N}\left\vert \lambda_{i}(G)-\lambda_{i}(H)\right\vert
,\qquad cs(G)=\min\{\sigma(G,H):H\not \cong G\},
\]
where the minimum is taken over all graphs $H$ of order $N$ \cite{AZ}.

The \emph{cocktail party graph} $CP_{n}=K_{2,2,\ldots,2}$ with $n$ parts is
the complete multipartite graph on $N=2n$ vertices whose parts all have size
two. It is $\left( 2n-2\right) $-regular, and its complement is the perfect
matching $nK_{2}.$ Its spectrum is
\begin{equation}
\operatorname{Spec}(CP_{n})=\{\,2n-2,\;0^{\,n},\;(-2)^{\,n-1}\,\}.
\label{eq_spec_CP}%
\end{equation}

Abdollahi and Zakeri \cite{AZ} proposed two conjectures on the cospectrality
of $CP_{n}.$ One is for the spectral distance $\sigma,$ and the other is for
the Euclidean distance $\sum_{i}(\lambda_{i}(G)-\lambda_{i}(H))^{2}.$ Let
\[
H_{n}:=(CP_{n-1}\nabla K_{1})\cdot K_{2}%
\]
be the graph obtained from the complete multipartite graph $K_{1,2,2,\ldots
,2}=CP_{n-1}\nabla K_{1}$ by attaching a pendant vertex to the vertex in the
part of size one.

\begin{conjecture}
[Manhattan, \cite{AZ}]\label{Conj_C1}For every $n\geq2,$ $cs(CP_{n}%
)=\sigma(CP_{n},H_{n}),$ and $H_{n}$ is the unique graph at which the
cospectrality of $CP_{n}$ is attained.
\end{conjecture}

The Euclidean conjecture of \cite{AZ} states that the Euclidean cospectrality
of $CP_{n}$ is attained exactly at the graph $CP_{n}-e,$ obtained by deleting
one edge. We do not consider it here. Besides
$H_{n},$ two further graphs play a role:
\[
B_{n}:=K_{1,3,2,\ldots,2}\ \ (n-2\text{ parts of size }2),\qquad H_{n}%
^{-}:=(CP_{n-1}\nabla K_{1})\cup K_{1}.
\]
Thus $B_{n}$ is $CP_{n}$ with one vertex moved from one part to another, and
$H_{n}^{-}$ is $H_{n}$ without its pendant edge. Note that $B_{2}=H_{2}.$ Our
main result is the following classification.

\begin{theorem}
[Main Theorem]\label{Tm_main}For every $n\geq2,$ up to isomorphism,
\[
\{H\not \cong CP_{n}:\left\vert H\right\vert =2n,\ \sigma(CP_{n},H)<2\}=\{H_{n},B_{n},H_{n}^{-}\}.
\]
\end{theorem}

Lemma \ref{Lem_C1_order} compares the three values $\sigma(CP_{n},\cdot)$ on
the right. Since $\sigma(CP_{n},H_{n})$ is the smallest of them, Theorem
\ref{Tm_main} implies Conjecture \ref{Conj_C1} for every $n\geq2.$ The proof of
Theorem \ref{Tm_main} is computer-assisted. The computations concern explicit
finite lists of graphs, and we describe them in Appendix \ref{App_comp}.

We work with the complement $G=\overline{H}$ of a competing graph $H.$ Section
\ref{Sec_prelim} collects the tools. These are an exact formula for
$\sigma(CP_{n},H),$ a lower bound $T$ in terms of the eigenvalues of
$J-A(G),$ and the hereditary class of admissible graphs. Section \ref{Sec_red}
proves Theorem \ref{Tm_main}. We compute the three graphs above and derive
necessary conditions. Then we show that, apart from components $K_{2},$ the
graph $G$ consists of one main component $G_{0}$ and a small forest $Z.$
Finally, we deduce Theorem \ref{Tm_main} from two theorems. The first treats
$G_{0}$ with a dominating vertex, and the second treats $G_{0}$ without one.
They are proved in Sections \ref{Sec_cone} and \ref{Sec_M}.

\section{Preliminaries}\label{Sec_prelim}

Throughout, $H$ is a graph on $2n$ vertices with eigenvalues $\mu_{1}%
\geq\cdots\geq\mu_{2n},$ $G:=\overline{H}$ and $\sigma(H):=\sigma(CP_{n},H).$
We write $J$ for the all-ones matrix and $\mathbf{1}$ for the all-ones vector.
For a symmetric matrix (or a graph, via its adjacency matrix), $n_{+}$ and
$n_{-}$ denote the numbers of positive and negative eigenvalues. For real $x$
let $x^{+}:=\max(x,0)$ and $x^{-}:=\max(-x,0).$ We put
\[
f(x):=\bigl\vert\left\vert x\right\vert -1\bigr\vert,\qquad g(x):=x+f(x),\qquad
\tau:=\tfrac{1+\sqrt{5}}{2}.
\]
The function $g$ is non-decreasing. It equals $-1$ on $(-\infty,-1],$ $1+2x$
on $[-1,0],$ $1$ on $[0,1]$ and $2x-1$ on $[1,\infty).$

We first give an exact formula for $\sigma(H).$ Then we bound $\sigma(H)$ from
below by a quantity $T$ which depends only on the eigenvalues of $J-A(G).$

\begin{lemma}
\label{Lem_C1_uv}Let $H$ be a graph on $2n$ vertices. Put $u(H)=\sum
_{i=2}^{n+1}\mu_{i}^{+}+\sum_{i=n+2}^{2n}(\mu_{i}+2)^{+}$ and $v(H)=\sum
_{i=2}^{n+1}\mu_{i}^{-}+\sum_{i=n+2}^{2n}(\mu_{i}+2)^{-}.$ Then
$\sigma(H)=2\max\bigl(u(H),v(H)\bigr)$ and $u(H)-v(H)=2n-2-\mu_{1}.$ In
particular, $\sigma(H)<2$ if and only if $u(H)<1$ and $v(H)<1.$
\end{lemma}

\begin{proof}
Let $\lambda_{i}$ be the eigenvalues of $CP_{n}$ and $e_{i}:=\lambda_{i}%
-\mu_{i}.$ Then $\sum_{i}e_{i}=\operatorname{tr}A(CP_{n})-\operatorname{tr}%
A(H)=0.$ By (\ref{eq_spec_CP}), $\sum_{i\geq2}e_{i}^{+}=v$ and $\sum_{i\geq
2}e_{i}^{-}=u.$ Hence $e_{1}=u-v$ and $\sigma=\sum_{i}\left\vert
e_{i}\right\vert =2\sum_{i}e_{i}^{+}=2(e_{1}^{+}+v)=2\max(u,v).$
\end{proof}

For a graph $W$ of order $N$ let $\xi_{1}(W)\geq\cdots\geq\xi_{N}(W)$ be the
eigenvalues of $M=J-A(W),$ and let $\delta(W):=N-1-\xi_{1}(W).$ For
$W=G=\overline{H}$ we have $J-A(G)=A(H)+I,$ so $\xi_{i}=\mu_{i}+1.$ Since
$\operatorname{tr}M=N,$ we have
\[
T(W):=\delta(W)+\sum_{i\geq2}f\bigl(\xi_{i}(W)\bigr)=-1+\sum_{i\geq
2}g\bigl(\xi_{i}(W)\bigr).
\]
If $W$ has order $2m$ and $\mu_{i}$ are the eigenvalues of $\overline{W},$
then $\sigma(CP_{m},\overline{W})\geq T(W).$ Indeed, for $i\geq2$ the $i$-th
eigenvalue of $CP_{m}$ lies in $\{0,-2\}.$ So its distance to $\mu_{i}$ is at
least $\operatorname{dist}(\mu_{i},\{0,-2\})=f(\mu_{i}+1)=f(\xi_{i}).$
Moreover, $\left\vert 2m-2-\mu_{1}\right\vert \geq2m-1-\xi_{1}.$

\begin{lemma}
[adding twin pairs]\label{Lem_P_mono}Let $G$ be a graph with at least one
vertex. Let $G^{\prime}$ be obtained from $G$ by adding two adjacent vertices
with the same neighbours in $G.$ Then $T(G^{\prime})\geq T(G).$
\end{lemma}

\begin{proof}
Let $p,p^{\prime}$ be the new vertices, $N=\left\vert G\right\vert ,$
$M=J-A(G)$ and $M^{\prime}=J-A(G^{\prime}).$ The rows $p$ and $p^{\prime}$ of
$M^{\prime}$ coincide outside $\{p,p^{\prime}\}.$ Moreover, $M_{pp}^{\prime}%
=M_{p^{\prime}p^{\prime}}^{\prime}=1$ and $M_{pp^{\prime}}^{\prime}=0.$ Hence
$e_{p}-e_{p^{\prime}}$ is an eigenvector of $M^{\prime}$ for $1.$ On its
orthogonal complement we use the basis $(e_{p}+e_{p^{\prime}})/\sqrt{2},$
$e_{v}$ $(v\in V(G)).$ In this basis $M^{\prime}$ is represented by a matrix
$M^{\prime\prime}$ which has $M$ as a principal submatrix. By interlacing, the
eigenvalues $\xi_{1}^{\prime\prime}\geq\cdots\geq\xi_{N+1}^{\prime\prime}$ of
$M^{\prime\prime}$ satisfy $\xi_{i}^{\prime\prime}\geq\xi_{i}(G)$ for $i\leq
N.$ As $\xi_{1}^{\prime\prime}\geq\xi_{1}(G)\geq M_{vv}=1,$ the eigenvalues
$\xi_{i}(G^{\prime}),$ $i\geq2,$ are $1,\xi_{2}^{\prime\prime},\ldots
,\xi_{N+1}^{\prime\prime}.$ Since $g$ is non-decreasing, $g(1)=1$ and
$g\geq-1,$ we get $T(G^{\prime})=-1+g(1)+\sum_{i=2}^{N+1}g(\xi_{i}%
^{\prime\prime})\geq\sum_{i=2}^{N}g\bigl(\xi_{i}(G)\bigr)-1=T(G).$
\end{proof}

In particular, Lemma \ref{Lem_P_mono} applies to components $K_{2}.$ Let $W$
be a graph of order $2m$ with at least one vertex. Then $\overline{W\cup
rK_{2}}=\overline{W}\nabla CP_{r}.$ We add the components $K_{2}$ one at a
time. The above bound and Lemma \ref{Lem_P_mono} then give, for all integers
$r\geq r_{0}\geq0,$
\begin{equation}
\sigma(CP_{m+r},\overline{W}\nabla CP_{r})\geq T(W\cup rK_{2})\geq
T(W\cup r_{0}K_{2}).\label{eq_K2_mono}%
\end{equation}

Finally, for a graph $W$ put
\[
\operatorname{cost}(W):=\max\Bigl(\sum_{i\geq2}\lambda_{i}(\overline{W}%
)^{+},\ \sum_{i}\bigl(\lambda_{i}(\overline{W})+2\bigr)^{-}\Bigr).
\]
We call $W$ \emph{admissible} if $\operatorname{cost}(W)<1,$ and we denote by
$\mathcal{P}$ the class of admissible graphs. In Lemma \ref{Lem_C1_nec}(v) we
will see that every induced subgraph of $G=\overline{H}$ is admissible if
$\sigma(H)<2.$ Let $\operatorname{def}(W):=\left\vert W\right\vert -\Delta
(W)$ be the \emph{deficiency} of $W.$ Thus $\operatorname{def}(W)=1$ if and
only if $W$ has a dominating vertex. A connected $W$ with
$\operatorname{def}(W)=d$ is a cone with $d-1$ missing spokes.

\begin{lemma}
\label{Lem_P_basic}\emph{(a)} The function $\operatorname{cost}$ does not
increase when we pass to an induced subgraph. In particular, $\mathcal{P}$ is
hereditary.

\emph{(b)} If $W\in\mathcal{P},$ then $\lambda_{2}(\overline{W})<1.$ Hence
$\overline{W}$ has no induced $2K_{2},$ and $W$ has no induced $C_{4}.$
\end{lemma}

\begin{proof}
(a) Let $F^{\prime}$ be an induced subgraph of $F$ on $k^{\prime}$ of its $k$
vertices. Interlacing gives $\lambda_{i}(F)\geq\lambda_{i}(F^{\prime}%
)\geq\lambda_{i+k-k^{\prime}}(F).$ Hence $\sum_{i\geq2}\lambda_{i}(F^{\prime
})^{+}\leq\sum_{i\geq2}\lambda_{i}(F)^{+}$ and $\sum_{i}(\lambda_{i}%
(F^{\prime})+2)^{-}\leq\sum_{i}(\lambda_{i}(F)+2)^{-}.$ We apply this to
$F=\overline{W}.$ (b) We have $\lambda_{2}(\overline{W})^{+}\leq
\operatorname{cost}(W)<1.$ On the other hand, $\lambda_{2}(2K_{2})=1,$ and
$\lambda_{2}$ does not increase when we pass to an induced subgraph.
\end{proof}

\section{Proof of the Main Theorem}\label{Sec_red}

\paragraph{The three graphs.}The complements of $H_{n},$ $B_{n}$ and
$H_{n}^{-}$ are $\overline{H_{n}}=K_{1}\cup(K_{1}\nabla(n-1)K_{2}),$
$\overline{B_{n}}=K_{3}\cup K_{1}\cup(n-2)K_{2}$ and $\overline{H_{n}^{-}%
}=K_{1}\nabla(K_{1}\cup(n-1)K_{2}).$ We compute their spectra and compare
their distances to $CP_{n}.$

\begin{proposition}
\label{Prop_C1_three}Let $n\geq3.$

\emph{(i)} $\operatorname{Spec}(H_{n})=\{x_{1},x_{2},0^{\,n-1},x_{3}%
,(-2)^{\,n-2}\},$ where $x_{1}>x_{2}>x_{3}$ are the roots of the cubic
$P(x)=x^{3}-(2n-4)x^{2}-(2n-1)x+(2n-4),$ with $x_{1}\in(2n-3,2n-2),$
$x_{2}\in(0,1)$ and $x_{3}\in(-2,0).$ Moreover, $\sigma(H_{n})=4n-4-2x_{1}<2.$

\emph{(ii)} $\operatorname{Spec}(B_{n})=\{\rho_{1},0^{\,n},\rho_{2}%
,(-2)^{\,n-3},\rho_{3}\},$ where $\rho_{1}>\rho_{2}>\rho_{3}$ are the roots of
$c_{n}(x)=x^{3}+(6-2n)x^{2}+(13-8n)x+(6-6n),$ with $\rho_{1}\in(0,2n-2),$
$\rho_{2}\in(-2,-1)$ and $\rho_{3}<-2.$ Moreover, $\sigma(B_{n})=4+2\rho
_{2}<2.$

\emph{(iii)} $\operatorname{Spec}(H_{n}^{-})=\{\kappa_{+},0^{\,n},\kappa
_{-},(-2)^{\,n-2}\}$ with $\kappa_{\pm}=n-2\pm\sqrt{n^{2}-2n+2}.$ Moreover,
$\sigma(H_{n}^{-})=2n-2\sqrt{n^{2}-2n+2}<2.$
\end{proposition}

\begin{proof}
Each of the three graphs has an induced $CP_{k}$ $(k=n-1,n-2,n-1)$ on a set
$X,$ and the vertices of $X$ have the same neighbours outside $X.$ If
$\{a,b\}$ is a part of $CP_{k},$ then $e_{a}-e_{b}$ is an eigenvector for $0.$
If $\{a,b\}$ and $\{a',b'\}$ are two parts, then
$e_{a}+e_{b}-e_{a'}-e_{b'}$ is an eigenvector for $-2.$ Consider the partition
of the vertex set into $X$ and the classes of vertices outside $X$ with the
same neighbourhood. This partition is equitable. The remaining eigenvalues are
those of its quotient matrix, and $0$ for each further vector $e_{a}-e_{b}$
with $a,b$ in one class.

(i) The classes outside $X$ are the apex $v$ and the pendant vertex $p.$ We
get $0^{\,n-1},$ $(-2)^{\,n-2}$ and the eigenvalues of
$\left(\begin{smallmatrix}2n-4 & 1 & 0\\ 2n-2 & 0 & 1\\ 0 & 1 & 0
\end{smallmatrix}\right),$ whose characteristic polynomial is $P.$ We have
$P(0)=2n-4>0,$ $P(1)=P(-2)=2-2n<0,$ $P(2n-3)=2-2n<0$ and $P(2n-2)>0.$ So the
roots are located as stated, and the sorted spectrum is $x_{1},x_{2}%
,0^{\,n-1},x_{3},(-2)^{\,n-2}.$ We compare it position by position with
(\ref{eq_spec_CP}). Since $x_{1}+x_{2}+x_{3}=2n-4,$ this gives
$\sigma=(2n-2-x_{1})+x_{2}+(x_{3}+2)=4n-4-2x_{1}.$ Finally, $x_{1}>2n-3$ gives
$\sigma<2.$

(ii) The classes outside $X$ are the parts of sizes $1$ and $3.$ We get
$0^{\,n},$ $(-2)^{\,n-3}$ and the eigenvalues of $\left(\begin{smallmatrix}0 &
3 & 2n-4\\ 1 & 0 & 2n-4\\ 1 & 3 & 2n-6\end{smallmatrix}\right),$ which are the
roots of $c_{n}.$ We have $c_{n}(-2)=2n-4>0,$ $c_{n}(-1)=-2,$ $c_{n}%
(0)=6-6n<0$ and $c_{n}(2n-2)=4(n-1)>0,$ and this locates the roots. Since
$\rho_{1}+\rho_{2}+\rho_{3}=2n-6,$ as in (i) we get
$\sigma=(2n-2-\rho_{1})+(\rho_{2}+2)+(-2-\rho_{3})=4+2\rho_{2}.$

(iii) The classes outside $X$ are the apex and the isolated vertex. We get
$0^{\,n},$ $(-2)^{\,n-2}$ and the roots $\kappa_{\pm}$ of
$x^{2}-(2n-4)x-(2n-2);$ one of the zeros comes from the isolated vertex. Since
$n-1<\sqrt{n^{2}-2n+2}<n,$ we have $\kappa_{+}<2n-2$ and $-2<\kappa_{-}<-1.$
Hence $\sigma=(2n-2-\kappa_{+})+(\kappa_{-}+2)=2n-2\sqrt{n^{2}-2n+2}<2.$
\end{proof}

For $n=2$ we have $H_{2}=B_{2}=K_{1,3},$ $\sigma(H_{2})=4-2\sqrt{3}$ and
$\sigma(H_{2}^{-})=4-2\sqrt{2}.$ The value $2$ is a natural barrier, since
$\sigma(CP_{n},CP_{n}-e)=2$ for every $n\geq2.$ This follows from the spectrum
of $CP_{n}-e.$ For $n=2$ it is $\{\pm\tau,\pm(\tau-1)\}.$ For $n\geq3$ it
consists of $0^{\,n-2},$ $(-2)^{\,n-3},$ $\tau-1,$ $-\tau$ and the roots of
$x^{3}-(2n-5)x^{2}-(6n-9)x-(2n-2),$ computed as in Proposition
\ref{Prop_C1_three}.

\begin{lemma}
\label{Lem_C1_order}For every $n\geq2,$ $\sigma(H_{n})<\sigma(H_{n}^{-}).$
For every $n\geq3,$ $\sigma(H_{n})<\sigma(B_{n})<\sigma(H_{n}^{-})$ $($recall
$B_{2}=H_{2}).$
\end{lemma}

\begin{proof}
For $n=2$ the first inequality reads $4-2\sqrt{3}<4-2\sqrt{2}.$ Let $n\geq3.$
By Proposition \ref{Prop_C1_three}, $\sigma(H_{n})<\sigma(H_{n}^{-})$ is
equivalent to $x_{1}=\lambda_{1}(H_{n})>\kappa_{+}=\lambda_{1}(H_{n}^{-}).$
This holds by the Perron--Frobenius theorem, since $H_{n}^{-}$ is a proper
spanning subgraph of the connected graph $H_{n}.$

For the second inequality put $L=2/(4n-5)\in(0,1).$ A direct computation
gives
\[
P(2n-3+L)=-\frac{2\bigl((n-2)(16n^{2}-24n+17)+5\bigr)}{(4n-5)^{3}}%
<0,\qquad c_{n}(-1-L)=\frac{2(4n-9)}{(4n-5)^{3}}>0.
\]
We have $P<0$ on $(x_{2},x_{1})$ and $P>0$ on $(x_{1},\infty).$ So the first
value gives $x_{1}>2n-3+L,$ i.e.\ $\sigma(H_{n})<2-2L.$ On $(-2,-1)$ the
polynomial $c_{n}$ changes sign exactly once, from positive to negative, at
$\rho_{2}.$ So the second value gives $\rho_{2}>-1-L,$ i.e.\ $\sigma
(B_{n})>2-2L.$ Finally, we reduce $c_{n}$ modulo $x^{2}-(2n-4)x-(2n-2),$ whose
roots are $\kappa_{\pm}.$ This gives $c_{n}(\kappa_{-})=(3-2n)\kappa_{-}%
+2-2n.$ This value is negative, because $\kappa_{-}>-1-1/(2n-3)$ by
$(n-1+1/(2n-3))^{2}>n^{2}-2n+2.$ Since $\kappa_{-}\in(-2,-1),$ we get
$\rho_{2}<\kappa_{-},$ i.e.\ $\sigma(B_{n})=4+2\rho_{2}<4+2\kappa_{-}%
=\sigma(H_{n}^{-}).$
\end{proof}

All three values are $2-1/n+O(n^{-2}).$ The differences in Lemma
\ref{Lem_C1_order} are of order $n^{-2},$ numerically close to $1/(4n^{2}).$
This is why Conjecture \ref{Conj_C1} is difficult.

\paragraph{Necessary conditions.}We now determine the components of $G$ for
an arbitrary graph $H$ with $\sigma(H)<2.$

\begin{lemma}
\label{Lem_C1_nec}Let $\sigma(H)<2,$ and let $G=\overline{H}$ have eigenvalues
$\nu_{1}\geq\cdots\geq\nu_{2n}.$ Then:

\emph{(i)} $2n-3<\mu_{1}<2n-1;$

\emph{(ii)} $\mu_{2}<1,$ $\mu_{n+1}>-1>\mu_{n+2}$ and $\mu_{2n}>-3;$

\emph{(iii)} $M=J-A(G)$ is nonsingular, $n_{+}(M)=n+1$ and $n_{-}(M)=n-1.$
Moreover, $A(G)$ has no null vector orthogonal to $\mathbf{1}.$ In particular,
$G$ has no \emph{false twins,} i.e.\ no two non-adjacent vertices with the
same neighbourhood;

\emph{(iv)} $\nu_{2}<2.$ Hence at most one component of $G$ has largest
eigenvalue $\geq2,$ and every other component is one of the Dynkin trees
$A_{k}=P_{k},$ $D_{k}$ $(k\geq4),$ $E_{6},E_{7},E_{8};$

\emph{(v)} every induced subgraph of $G,$ in particular every component of
$G,$ is admissible.
\end{lemma}

\begin{proof}
By Lemma \ref{Lem_C1_uv}, $u,v<1.$ (i) $\left\vert 2n-2-\mu_{1}\right\vert
=\left\vert u-v\right\vert <1.$ (ii) We have $\mu_{2}^{+}\leq u,$ $\mu
_{n+1}^{-}\leq v,$ $(\mu_{n+2}+2)^{+}\leq u$ and $(\mu_{2n}+2)^{-}\leq v.$
(iii) The eigenvalues of $M=A(H)+I$ are $\xi_{i}=\mu_{i}+1.$ By (i) and (ii),
$\xi_{1}>0$ and $\xi_{n+1}>0>\xi_{n+2}.$ If $A(G)x=\lambda x$ and
$\mathbf{1}^{\top}x=0,$ then $Mx=Jx-A(G)x=-\lambda x.$ For $\lambda=0$ this is
impossible, as $M$ is nonsingular. False twins $a,b$ would give such a null
vector $e_{a}-e_{b}.$

(iv) We have $A(G)=(-I-A(H))+J.$ Hence $\nu_{2}\leq\lambda_{1}(-I-A(H))=-1-\mu
_{2n}<2.$ Two components with largest eigenvalue $\geq2$ would give $\nu
_{2}\geq2.$ The connected graphs with largest eigenvalue $<2$ are exactly the
Dynkin trees (Smith's theorem, \cite{CvetRowSim}).

(v) By Lemma \ref{Lem_P_basic}(a), it suffices to show
$\operatorname{cost}(G)<1.$ By (ii), $\sum_{i\geq2}\mu_{i}^{+}=\sum
_{i=2}^{n+1}\mu_{i}^{+}\leq u<1$ and $\sum_{i}(\mu_{i}+2)^{-}=\sum_{i\geq
n+2}(\mu_{i}+2)^{-}\leq v<1.$
\end{proof}

\begin{lemma}
\label{Lem_C1_struct}Let $\sigma(H)<2$ and $H\not \cong CP_{n}.$ Write
$G=\overline{H}=W\cup rK_{2}$ $(r\geq0),$ where $W$ has no component $K_{2}.$
Then $T(W)<2$ and $W$ has no false twins. Moreover, $W=G_{0}\cup Z,$ where
the \emph{main component} $G_{0}$ is connected and admissible with
$\lambda_{1}(G_{0})\geq2.$ Here $Z=Z_{1}\cup Z_{2}$ with
$Z_{1}\in\{\emptyset,P_{6},E_{6},P_{6}\cup E_{6}\}$ and
$Z_{2}\in\{\emptyset,K_{1},P_{5}\}$ (twelve graphs $Z$).
\end{lemma}

\begin{proof}
As $H\not \cong CP_{n},$ we have $W\neq\emptyset.$ By (\ref{eq_K2_mono}) with
$r_{0}=0,$ $T(W)\leq\sigma(H)<2.$ False twins of $W$ would be false twins of
$G,$ and $G$ has none by Lemma \ref{Lem_C1_nec}(iii). By Lemma
\ref{Lem_C1_nec}(iv),(v), all components of $G$ except at most one are
admissible Dynkin trees. Such a component of $W$ is not $P_{2}=K_{2}.$ It is
also not $P_{3}$ or $D_{k}$ $(k\geq4),$ since these have two leaves at a
common vertex, i.e.\ false twins. We compute $\operatorname{cost}%
(E_{7})=1.199\ldots$ and $\operatorname{cost}(P_{7})=1.252\ldots.$ Every
longer path contains $P_{7}$ as an induced subgraph, and $E_{8}$ contains
$E_{7}.$ As $\mathcal{P}$ is hereditary, the remaining candidates are
$P_{1},P_{4},P_{5},P_{6},E_{6}.$

Let $C$ be such a component of $W.$ Let $x$ be an eigenvector of $A(C)$ with
eigenvalue $\lambda$ and $\mathbf{1}^{\top}x=0,$ extended by zero to $V(G).$
As in the proof of Lemma \ref{Lem_C1_nec}(iii), $Mx=-\lambda x$ for
$M=J-A(G).$ The same holds for the vector equal to $y$ on $C$ and to $-y$ on
$C^{\prime},$ if $C^{\prime}\cong C$ is another component and $y$ is any
eigenvector of $A(C)$ with eigenvalue $\lambda.$ We choose such
vectors orthonormally. This gives a family of eigenvalues of $M,$ which we
call \emph{fixed}. They occur in $\operatorname{Spec}M$ with their
multiplicities. They have absolute value $<2,$ while $\xi_{1}=\mu_{1}+1>2n-2\geq2$
by Lemma \ref{Lem_C1_nec}(i). By Lemma \ref{Lem_C1_nec}(iii), $0$ is not
fixed. A positive fixed eigenvalue is some $\xi_{i}$ with $2\leq i\leq n+1,$
and a negative one is some $\xi_{i}$ with $i\geq n+2.$ We write the quantity
$u(H)$ of Lemma \ref{Lem_C1_uv} in terms of $\xi_{i}=\mu_{i}+1$ and get
\[
1>u(H)\geq\sum_{\xi\text{ fixed}}\theta(\xi),\qquad\theta(\xi):=(\xi
-1)^{+}\ (\xi>0),\quad\theta(\xi):=(1+\xi)^{+}\ (\xi<0).
\]
\emph{$P_{4}$.} The vectors $(1,\tau-1,1-\tau,-1)$ and $(1,-\tau,\tau,-1)$
are eigenvectors of $A(P_{4})$ for $\tau-1$ and $-\tau,$ with coordinate sum
$0.$ So $1-\tau$ and $\tau$ are fixed, and $u(H)\geq(2-\tau)+(\tau-1)=1,$ a
contradiction.

\emph{$P_{6},E_{6}$.} Suppose that two components are isomorphic to $C\in
\{P_{6},E_{6}\}.$ Then the fixed eigenvalues include $-\lambda$ for every
$\lambda\in\operatorname{Spec}C.$ As $\operatorname{Spec}C$ is symmetric,
$u(H)\geq\sum_{\lambda\in\operatorname{Spec}C,\,\lambda>0}\left\vert
\lambda-1\right\vert .$ The positive eigenvalues are $2\cos(j\pi/7),$
$j=1,2,3,$ for $P_{6},$ and $2\cos(\pi/12),1,2\cos(5\pi/12)$ for $E_{6}.$ The
sum is $1.603\ldots$ resp. $2\cos(\pi/12)-2\cos(5\pi/12)=\sqrt{2},$ a
contradiction.

\emph{$K_{1},P_{5}$.} Each of them has a null vector with non-zero coordinate
sum, namely $(1)$ resp. $(1,0,-1,0,1).$ Suppose that there are two such
components. Then a combination of these vectors is a null vector of $A(G)$
orthogonal to $\mathbf{1},$ which contradicts Lemma \ref{Lem_C1_nec}(iii). So
at most one of them occurs.

Finally, suppose that no component of $W$ has $\lambda_{1}\geq2.$ Then
$W=Z_{1}\cup Z_{2}$ has even order, so $Z_{2}=\emptyset$ and $W\in
\{P_{6},E_{6},P_{6}\cup E_{6}\}.$ A direct computation gives
$T(W)=2.826\ldots,$ $2.661\ldots,$ $5.806\ldots,$ respectively, a
contradiction.
\end{proof}

\paragraph{The main component.}Lemma \ref{Lem_C1_struct} leaves two cases.
Either $\operatorname{def}(G_{0})=1,$ i.e.\ $G_{0}=K_{1}\nabla Y$ for a graph
$Y,$ or $\operatorname{def}(G_{0})\geq2.$ The following two theorems treat
these cases. They are proved in Sections \ref{Sec_cone} and \ref{Sec_M},
respectively.

\begin{theorem}
[cones]\label{Tm_cone}Let $m\geq2,$ let $Y$ be a graph and let $Z$ be one of
the twelve graphs of Lemma \ref{Lem_C1_struct}. Let $G=(K_{1}\nabla Y)\cup Z$
be a graph of order $2m$ without false twins. If $\overline{G}\notin
\{H_{m},H_{m}^{-}\},$ then $T(G)\geq2.$
\end{theorem}

\begin{theorem}
[missing spokes]\label{Tm_M}Let $G=G_{0}\cup Z$ be a graph of order $2m$
without false twins. Here $G_{0}$ is a connected admissible graph with
$\lambda_{1}(G_{0})\geq2$ and $\operatorname{def}(G_{0})\geq2,$ and $Z$ is
one of the twelve graphs of Lemma \ref{Lem_C1_struct}. Then $T(G)\geq2.$
\end{theorem}

\begin{proof}
[Proof of Theorem \ref{Tm_main}]Let $\sigma(H)<2$ and $H\not \cong CP_{n}.$
Let $G=\overline{H}=W\cup rK_{2},$ $W=G_{0}\cup Z$ and $\left\vert
W\right\vert =2m$ be as in Lemma \ref{Lem_C1_struct}. Then $m\geq2,$ as
$\lambda_{1}(G_{0})\geq2.$ Since $T(W)<2,$ Theorem \ref{Tm_M} shows that
$\operatorname{def}(G_{0})=1,$ i.e.\ $G_{0}=K_{1}\nabla Y.$ Then Theorem
\ref{Tm_cone} gives $\overline{W}\in\{H_{m},H_{m}^{-}\}.$ If $r=0,$ then
$H=\overline{W}\in\{H_{n},H_{n}^{-}\}.$ Let $r\geq1.$ Then $H=\overline
{W}\nabla CP_{r},$ and $\sigma(H)\geq T(W\cup K_{2})$ by (\ref{eq_K2_mono}).

First let $\overline{W}=H_{m}^{-},$ i.e.\ $W=K_{1}\nabla(K_{1}\cup
(m-1)K_{2}),$ or $\overline{W}=H_{m}$ with $m\geq3,$ i.e.\ $W=K_{1}\cup
(K_{1}\nabla(m-1)K_{2}).$ We delete $m-2$ resp. $m-3$ components $K_{2}$ of
$Y.$ Their vertices have the apex as their only other neighbour. Lemma
\ref{Lem_P_mono} and a direct computation give $T(W\cup K_{2})\geq
T(K_{2}\cup(K_{1}\nabla(K_{1}\cup K_{2})))=2.021\ldots$ resp. $T(W\cup
K_{2})\geq T(K_{1}\cup K_{2}\cup(K_{1}\nabla2K_{2}))=2.356\ldots.$ So
$\sigma(H)>2.$ Otherwise $m=2$ and $\overline{W}=H_{2},$ i.e.\ $W=K_{3}\cup
K_{1}.$ Then $H=H_{2}\nabla CP_{n-2}=B_{n}.$ Conversely, $H_{n},B_{n},H_{n}^{-}$ satisfy $\sigma<2.$ For $n\geq3$ this is
Proposition \ref{Prop_C1_three}, and for $n=2$ it follows from the values
given after it.
\end{proof}

\section{Cones}\label{Sec_cone}

In this section we prove Theorem \ref{Tm_cone}. Let $G=(K_{1}\nabla Y)\cup Z$
be as there, with apex $o.$ Let $M=J-A(G),$ $\xi_{i}=\xi_{i}(G)$ and
$\delta=\delta(G).$ The hard case is $Z\in\{\emptyset,K_{1}\},$ since it
contains $\overline{H_{m}}$ and $\overline{H_{m}^{-}}.$ Lemma
\ref{Lem_cone_delta} prepares this case. The other ten graphs $Z$ are treated
by a perturbation argument, based on the functional $\eta_{1}$ below.

\begin{lemma}
\label{Lem_cone_delta}Let $Z\in\{\emptyset,K_{1}\},$ and let $G=(K_{1}\nabla
Y)\cup Z$ be a graph of even order without false twins. Assume that $Y$ has no
component $K_{2}$ and $\left\vert Y\right\vert \geq3.$ If $Z=K_{1},$ assume
also that $Y\notin\{P_{4},K_{1}\cup K_{3},K_{1}\cup P_{5}\}.$ Then
$\delta\geq1.$
\end{lemma}

\begin{proof}
Let $N=\left\vert G\right\vert ,$ $q=\left\vert Y\right\vert $ and
$\xi=\xi_{1}.$ Let $d_{v}$ be the degrees in $Y.$ Let $i_{0}$ and $k$ be the
numbers of vertices of $Y$ with $d_{v}=0$ and with $d_{v}\geq2,$
respectively. Two isolated vertices of $Y$ would be false twins of $G,$ and so
would two leaves of $Y$ with a common neighbour. Two adjacent leaves would
form a component $K_{2}.$ So $i_{0}\leq1,$ and the leaves of $Y$ have
distinct neighbours, of degree at least $2.$ Hence $q\leq i_{0}+2k.$

Suppose that $\xi>N-2=q-1+\left\vert Z\right\vert .$ First let $Z=\emptyset.$
Then $M=(1)\oplus(J-A(Y))$ has the row sums $1$ and $q-d_{v}.$ These are at
most $N-2$ unless $d_{v}=0.$ As $\xi$ is at most the largest row sum of the
non-negative matrix $M,$ we get $i_{0}=1.$ So $q\neq3,$ since otherwise the
other two vertices would form a component $K_{2}$ or be isolated. Since
$q=N-1$ is odd, $q\geq5.$ In this case we take an eigenvector $p\geq0$ of $M$
for $\xi$ supported on $V(Y).$ We put $p_{z}:=0$ and $s:=\mathbf{1}^{\top
}p_{Y}>0.$ Now let $Z=K_{1}=\{z\},$ and let $p\geq0,$ $p\neq0,$ be an
eigenvector of $M$ for $\xi.$ The rows $o$ and $z$ of $Mp=\xi p$ read $\xi
p_{o}=p_{o}+p_{z}$ and $\xi p_{z}=p_{o}+s,$ with $s:=p_{z}+\mathbf{1}^{\top
}p_{Y}.$ Hence $p_{z}=(\xi-1)p_{o}$ and $s=(\xi^{2}-\xi-1)p_{o}>0.$ (If
$p_{o}=0,$ then $p_{z}=s=0$ and $p=0.$) In both cases the rows of $V(Y)$ give
$(\xi I+A(Y))x=\mathbf{1}$ for $x:=p_{Y}/s\geq0,$ where $\mathbf{1}^{\top
}x=1-p_{z}/s.$ We sum the coordinates of this equation:
\[
\sum_{v}(d_{v}-1)x_{v}=q-(\xi+1)\mathbf{1}^{\top}x=q-\xi-1+\frac{(\xi
+1)p_{z}}{s}<\frac{\left\vert Z\right\vert (\xi+2)}{\xi^{2}}.
\]
Indeed, if $Z=\emptyset,$ then $q-\xi-1<0.$ If $Z=K_{1},$ then $q-\xi<0$ and
$(\xi+1)p_{z}/s=(\xi^{2}-1)/(\xi^{2}-\xi-1)=1+\xi/(\xi^{2}-\xi-1)<1+(\xi
+2)/\xi^{2},$ as $\xi>q\geq4.$ On the other hand, $A(Y)x\geq0$ gives
$x\leq\mathbf{1}/\xi.$ Hence $(A(Y)x)_{v}\leq d_{v}/\xi$ and $x_{v}\geq
(\xi-d_{v})/\xi^{2},$ with equality if $d_{v}=0.$ Therefore
\begin{equation}
\textstyle\sum_{d_{v}\geq2}(d_{v}-1)(\xi-d_{v})<i_{0}\xi+\left\vert
Z\right\vert (\xi+2).\label{eq_cone_count}%
\end{equation}
All $d_{v}$ are at most $\xi-1;$ if $Z=\emptyset,$ then $d_{v}\leq q-2$ as
$i_{0}=1.$ Moreover, $(d-1)(\xi-d)\geq\xi-2$ for $2\leq d\leq\xi-1.$ So
$k(\xi-2)<i_{0}\xi+\left\vert Z\right\vert (\xi+2).$ If $Z=\emptyset,$ this
gives $k<\xi/(\xi-2)<2,$ as $\xi>q-1\geq4.$ Then $q\leq i_{0}+2k\leq3,$ a
contradiction.

Let $Z=K_{1}.$ Then $q\geq4$ is even, and $k<i_{0}+1+(2i_{0}+4)/(\xi-2)$ with
$\xi>q.$ Let $q\geq6.$ Then $k\leq1$ if $i_{0}=0,$ and $k\leq3$ if $i_{0}=1.$
As $q\leq i_{0}+2k,$ we get $q=6,$ $i_{0}=1$ and $k=3.$ So $Y=K_{1}\cup
Y_{1},$ where $Y_{1}$ has two leaves with distinct neighbours $a,b$ and a fifth
vertex. This vertex has degree at least $2,$ so it is adjacent to $a$ and
$b.$ If $a\not\sim b,$ then $Y=K_{1}\cup P_{5}.$ If $a\sim b,$ the left side
of (\ref{eq_cone_count}) is $5\xi-14>2\xi+2,$ a contradiction. Finally let
$q=4.$ If $i_{0}=0,$ then $k<1+4/(\xi-2)<3$ and $q\leq2k$ give $k=2.$ The two
leaves are attached to the other two vertices, so these two vertices are
adjacent, and $Y=P_{4}.$ If $i_{0}=1,$ the other three vertices of $Y$ induce
$P_{3}$ or $K_{3}.$ The graph $P_{3}$ has two leaves with a common neighbour,
so $Y=K_{1}\cup K_{3}.$
\end{proof}

For a graph $W$ let $\eta_{1}(W):=\sum_{j}f\bigl(\lambda_{j}(W)\bigr).$ It is additive over
components. Moreover, $\eta_{1}(Y)\geq1$ for every non-empty graph $Y$ without
component $K_{2}.$ By additivity, it suffices to show $\eta_{1}(W)\geq1$ for
every connected $W\not \cong K_{2}.$ For $W=K_{w}$ we have
$\eta_{1}(W)=\left\vert w-2\right\vert \geq1.$ Otherwise $W$ contains an
induced $P_{3},$ so $\lambda_{\min}(W)\leq\lambda_{\min}(P_{3})=-\sqrt{2}$ by
interlacing. If $W=P_{3},$ then $\eta_{1}(W)=2\sqrt{2}-1.$ If $\left\vert
W\right\vert \geq4,$ a spanning tree of $W$ contains $P_{4}$ or $K_{1,3}.$ So
$\lambda_{1}(W)\geq\min(\tau,\sqrt{3})=\tau,$ and $\eta_{1}(W)\geq
(\tau-1)+(\sqrt{2}-1)>1.$

\begin{proof}
[Proof of Theorem \ref{Tm_cone}]Write $Y=Y_{0}\cup aK_{2},$ where $Y_{0}$ has
no component $K_{2}.$ Let $G^{\prime}=(K_{1}\nabla Y_{0})\cup Z,$ of order
$2m^{\prime},$ and let $M^{\prime}=J-A(G^{\prime})$ and $\xi_{i}^{\prime}%
=\xi_{i}(G^{\prime}).$ The vertices of a component $K_{2}$ of $Y$ have the
neighbourhood $\{o\}$ in $G-V(K_{2}).$ So Lemma \ref{Lem_P_mono} gives
$T(G)\geq T(G^{\prime}),$ and it suffices to show $T(G^{\prime})\geq2.$

\emph{Case $Z\in\{\emptyset,K_{1}\}.$} Here
$\overline{G}\notin\{H_{m},H_{m}^{-}\}$ means $Y\not \cong(m-1)K_{2}$ for
$Z=K_{1}$ and $Y\not \cong K_{1}\cup(m-1)K_{2}$ for $Z=\emptyset.$ By this
and by parity, $Y_{0}\notin\{\emptyset,K_{1}\}.$ The deleted vertices are adjacent only to $o$ and
to each other. So every vertex of $G^{\prime}$ other than $o$ has the same
neighbourhood in $G^{\prime}$ as in $G.$ The only possible non-neighbour $z$
of $o$ has $N_{G^{\prime}}(z)=\emptyset.$ Hence $G^{\prime}$ has no false
twins. In particular, $Y_{0}\neq2K_{1},$ so $\left\vert Y_{0}\right\vert
\geq3.$ If $Z=K_{1}$ and $Y_{0}\in\{P_{4},K_{1}\cup K_{3},K_{1}\cup P_{5}\},$
a direct computation gives $T(G^{\prime})=3.008\ldots,$ $2.042\ldots,$
$3.990\ldots,$ respectively. Otherwise Lemma \ref{Lem_cone_delta} gives
$\delta(G^{\prime})\geq1.$

Now $T(G^{\prime})=-1+\sum
_{i\geq2}g(\xi_{i}^{\prime}),$ $\xi_{1}^{\prime}\geq1$ and $\sum_{i\geq2}%
\xi_{i}^{\prime}=2m^{\prime}-\xi_{1}^{\prime}=\delta(G^{\prime})+1.$ First
let $n_{-}(M^{\prime})\leq m^{\prime}-2.$ We use $g(x)\geq-1+2[x\geq0]$ and
$\#\{i\geq2:\xi_{i}^{\prime}\geq0\}=2m^{\prime}-1-n_{-}(M^{\prime}).$ This
gives $T(G^{\prime})\geq2(m^{\prime}-1-n_{-}(M^{\prime}))\geq2.$ Otherwise
$n_{+}(M^{\prime})\leq2m^{\prime}-n_{-}(M^{\prime})\leq m^{\prime}+1.$ Now we
use $g(x)\geq2x-1+2[x\leq0]$ and $\#\{i\geq2:\xi_{i}^{\prime}\leq
0\}=2m^{\prime}-n_{+}(M^{\prime}).$ This gives $T(G^{\prime})\geq
2\delta(G^{\prime})+2(m^{\prime}+1-n_{+}(M^{\prime}))\geq2.$

\emph{Case $Z\notin\{\emptyset,K_{1}\}.$} Now $Z$ is one of the ten remaining
graphs of Lemma \ref{Lem_C1_struct}, and $Y_{0}$ may be empty. In the basis
$(V(Y_{0}),o,V(Z))$ one checks entrywise that
\[
M^{\prime}=D+hh^{\top},\qquad D=\bigl(-A(Y_{0})\bigr)\oplus D^{\prime}(Z),\qquad
D^{\prime}(Z)=%
\begin{pmatrix}
1 & \mathbf{1}^{\top}\\
\mathbf{1} & -A(Z)
\end{pmatrix}
,
\]
with $h=(\mathbf{1}_{Y_{0}},0,\mathbf{1}_{Z})$ and $\left\Vert h\right\Vert
^{2}=2m^{\prime}-1.$ Put $c_{0}(Z):=\sum_{d\in\operatorname{Spec}D^{\prime}(Z)}f(d),$
$d_{1}^{\prime}(Z):=\lambda_{\max}(D^{\prime}(Z))\geq1$ and
$\operatorname{base}(Z):=c_{0}(Z)-g(d_{1}^{\prime}(Z)).$ Let $d_{1}\geq
d_{2}\geq\cdots$ be the eigenvalues of $D.$ Then $d_{1}=\max(d_{1}^{\prime
}(Z),\beta),$ where $\beta:=-\lambda_{\min}(Y_{0}),$ and $\beta:=0$ if
$Y_{0}=\emptyset.$ As $\operatorname{tr}M^{\prime}=\operatorname{tr}%
D+\left\Vert h\right\Vert ^{2},$ we have $T(G^{\prime})=\sum_{i\geq2}%
g(\xi_{i}^{\prime})-\sum_{i}d_{i}.$ By Weyl's inequality, $\xi_{i}^{\prime
}\geq d_{i},$ since $hh^{\top}\succeq0.$ Moreover, $g$ is non-decreasing,
$g(d)-d=f(d)$ and $f$ is even. Hence
\[
T(G^{\prime})\geq\sum_{i\geq2}f(d_{i})-d_{1}=\eta_{1}(Y_{0})+c_{0}%
(Z)-g(d_{1}).
\]
If $\beta>d_{1}^{\prime}(Z)\geq1,$ then $\lambda_{1}(Y_{0})\geq\beta$ by the
Perron--Frobenius theorem. So $\eta_{1}(Y_{0})\geq2\beta-2=g(\beta)-1$ and
$T(G^{\prime})\geq c_{0}(Z)-1.$ Otherwise $T(G^{\prime})\geq
\operatorname{base}(Z)+\eta_{1}(Y_{0}).$ We computed the spectra of the ten matrices $D^{\prime}(Z).$ They
give $c_{0}(Z)>4.8$ for all of them. Further, $\operatorname{base}%
(Z)=1.236\ldots,$ $1.981\ldots,$ $1.729\ldots$ for $Z=P_{5},P_{6},E_{6},$ and
$\operatorname{base}(Z)\geq2.232$ for the other seven, with the minimum at
$E_{6}\cup K_{1}.$ As $\eta_{1}(Y_{0})\geq1$ for $Y_{0}\neq\emptyset,$ we get
$T(G^{\prime})\geq2$ unless $Y_{0}=\emptyset$ and $\operatorname{base}(Z)<2.$
In that case $G^{\prime}=K_{1}\cup Z$ has even order, so $\left\vert
Z\right\vert $ is odd, $Z=P_{5}$ and $T(G^{\prime})=T(K_{1}\cup
P_{5})=2.443\ldots$
\end{proof}

\section{Missing spokes}\label{Sec_M}

In this section we prove Theorem \ref{Tm_M}. Let $G=G_{0}\cup Z$ be as there,
and let $q:=\left\vert G_{0}\right\vert .$ Using a census of small admissible
graphs, we first show that $\operatorname{def}(G_{0})=2$ if $q\geq15.$ Then we
show that $G_{0}$ has a very restricted form (Lemma \ref{Lem_M_Ftype}). This
allows us to delete pairs of twin vertices until $q\leq14.$

It is easy to generate the class $\mathcal{P}$ by computer. Since
$\mathcal{P}$ is hereditary, every member of order $q$ is a one-vertex
extension of a member of order $q-1.$ A branch can be cut as soon as the cost
of an induced subgraph reaches $1.$ We generated the set $\mathcal{P}_{q}$ of
members of order $q$ for $q\leq15$ (Appendix \ref{App_lists}). We record:

\begin{description}
\item[(Cens15)] every connected member of $\mathcal{P}_{15}$ has deficiency at
most $2.$
\end{description}

\noindent This was also verified by a separate exhaustive computation, which
does not use the generation of $\mathcal{P}_{q}$ (Appendix \ref{App_lists}).

\begin{lemma}
\label{Lem_C1_RL}Let $G$ be a connected graph on $q\geq10$ vertices without
induced $C_{4},$ with $\operatorname{def}(G)\geq3.$ Then $G$ has a vertex $u$
such that $G-u$ is connected and $\operatorname{def}(G-u)\geq3.$
\end{lemma}

\begin{proof}
If $\operatorname{def}(G)\geq4,$ any non-cut vertex $u$ works. Let
$\operatorname{def}(G)=3$ and $F:=\overline{G}.$ Let $S_{2}$ be the set of
vertices of $F$-degree $2,$ let $B:=N_{F}(S_{2}),$ and let $C$ be the set of
cut vertices of $G.$ A vertex $u\notin B\cup C$ works. Indeed, $G-u$ is
connected. Since $u$ has no $F$-neighbour in $S_{2},$ every vertex of $G-u$
has at least two non-neighbours in $G-u.$ Suppose that $V(G)=B\cup C.$

First, $\left\vert C\right\vert \leq3.$ Fix $w\in S_{2}$ with $N_{F}%
(w)=\{x,y\}.$ For $c\in C\setminus\{w\},$ the graph $G-c$ has a component
$K\not\ni w.$ So $K\subseteq\{x,y\},$ and $c$ is the only $G$-neighbour of $K$
outside $K.$ At most two of $\{x\},\{y\},\{x,y\}$ occur as such $K.$ Indeed,
$\{x,y\}$ occurs only if $x\sim_{G}y.$ In that case $\{x\}$ and $\{y\}$ do not
both occur, since otherwise $\{x,y\}$ would be a component of $G.$

The graph $F$ has no induced $2K_{2},$ and $\Delta(F[S_{2}])\leq2.$ So
$F[S_{2}]$ has at most one component $Q$ with an edge. It is a path on at
most $4$ vertices, with at most two $F$-neighbours outside $Q,$ or a cycle on
at most $5$ vertices. So $\left\vert N_{F}(Q)\right\vert \leq6.$ The vertices
of $S_{2}^{\ast}:=S_{2}\setminus Q$ have no $F$-neighbour in $S_{2}.$ So
$S_{2}^{\ast}\cap B=\emptyset$ and $S_{2}^{\ast}\subseteq C.$ We have
$q-3\leq\left\vert B\right\vert \leq2\left\vert S_{2}\right\vert $ and
$q\geq10,$ so $\left\vert S_{2}\right\vert \geq4>\left\vert C\right\vert .$
Hence $Q\neq\emptyset.$ Moreover, $S_{2}^{\ast}\neq\emptyset,$ since
otherwise $q-3\leq\left\vert N_{F}(Q)\right\vert \leq6.$

Let $c\in S_{2}^{\ast}.$ For $s\in S_{2}\setminus\{c\},$ the components of
$G-c$ not containing $s$ consist of $F$-neighbours of $s.$ So they have at
most two vertices in total. Hence, as $q-1>4,$ one component contains
$S_{2}\setminus\{c\}.$ The set $R\neq\emptyset$ of vertices in the other
components lies in $N_{F}(s)$ for every $s\in S_{2}\setminus\{c\}.$ A vertex
of $R$ has all its $G$-neighbours in $R\cup\{c\}.$ So its $F$-degree is at
least $q-3>2,$ and $R\cap S_{2}=\emptyset.$ Fix $s\in Q.$ It has an
$F$-neighbour in $Q,$ so $R=N_{F}(s)\setminus S_{2}=\{r\}$ for a vertex $r$
which does not depend on $c.$ Moreover, $N_{G}(r)=\{c\}.$ Thus $S_{2}^{\ast
}=\{c\}.$ Every vertex of $Q$ is $F$-adjacent to $r\notin S_{2},$ so it has
exactly one $F$-neighbour in $Q.$ So $F[Q]=K_{2}$ and $\left\vert
S_{2}\right\vert =3,$ a contradiction.
\end{proof}

If $q\geq15,$ we will see that $G_{0}$ has exactly one missing spoke. We now
describe such graphs. Let $\operatorname{def}(G_{0})=2.$ Let $o$ be a vertex
of $G_{0}$ of degree $q-2,$ and let $w$ be its non-neighbour. Let
$Y:=G_{0}-o-w$ and $S:=N_{G_{0}}(w)\subseteq V(Y).$ Then $S\neq\emptyset.$
If $G_{0}\in\mathcal{P},$ then $S$ is a clique. Indeed, two non-adjacent
$s,s^{\prime}\in S$ would give an induced $C_{4}$ $o,s,w,s^{\prime}$ (Lemma
\ref{Lem_P_basic}(b)). We say that $G_{0}$ is \emph{F-type with respect to}
$o$ if
\[
Y\cong aK_{1}\cup bK_{2}\cup tK_{3}\qquad(a,b\geq0,\ t\in\{0,1\}).
\]
From the lists of Appendix \ref{App_lists} we read off:

\begin{description}
\item[(Cens14)] Every connected member of $\mathcal{P}_{14}$ with deficiency
$2$ is F-type with respect to every vertex of degree $12.$
\end{description}

\begin{lemma}
\label{Lem_M_Ftype}Assume \emph{(Cens14)}. Let $G_{0}\in\mathcal{P}$ be
connected with $\operatorname{def}(G_{0})=2$ and $q\geq15.$ Then $G_{0}$ is
F-type with respect to every vertex $o$ of degree $q-2.$
\end{lemma}

\begin{proof}
Suppose not. Then there is a \emph{witness} $U\subseteq V(Y)$: the vertices of
an induced $P_{3}$ or $K_{4}$ of $Y,$ or of two components $K_{3}$ of $Y.$ Fix
$s\in S.$ Every $x\in S$ is adjacent to $o$ and $w$ and has degree at most
$q-2.$ So $x$ has a non-neighbour $\bar{x}\in V(Y),$ and $\bar{x}\notin S$ as
$S$ is a clique. We build a set $Y^{\prime}\subseteq V(Y)$ in three steps.
First, $Y^{\prime}$ contains $U\cup\{s\}.$ Second, it contains $\bar{x}$ for
every $x\in S\cap(U\cup\{s\}).$ So far $Y^{\prime}$ has at most $11$
vertices, since $S\cap U$ is a clique of $Y[U].$ Third, we add further
vertices, first from $V(Y)\setminus S$ and then from $S,$ until $\left\vert
Y^{\prime}\right\vert =12;$ note that $\left\vert Y\right\vert \geq13.$ Then
$\bar{x}\in Y^{\prime}$ for every $x\in S\cap Y^{\prime}.$ Indeed, if $x$ was
added in the third step, then $V(Y)\setminus S\subseteq Y^{\prime}.$

Let $G^{\prime}:=G_{0}[\{o,w\}\cup Y^{\prime}]\in\mathcal{P}.$ It is
connected, since $o\sim Y^{\prime}\ni s\sim w.$ The vertex $o$ has degree
$12.$ No vertex is dominating: $o\not\sim w,$ $Y^{\prime}\setminus
S\not\sim w,$ and the vertices of $S\cap Y^{\prime}$ are not dominating by the
choice of $Y^{\prime}.$ So $\operatorname{def}(G^{\prime})=2.$ As
$G^{\prime}-o-w=Y[Y^{\prime}],$ the witness shows that $G^{\prime}$ is not
F-type with respect to $o.$ (Two components $K_{3}$ of $Y$ inside $Y^{\prime}$
are components of $Y[Y^{\prime}].$) This contradicts (Cens14).
\end{proof}

\begin{proof}
[Proof of Theorem \ref{Tm_M}]If $q\leq14,$ then $G$ belongs to an explicit
finite list. For this list $T(G)\geq2$ was checked by computer (Appendix
\ref{App_small}). Let $q\geq15.$ Then $\operatorname{def}(G_{0})=2.$
Otherwise, $G_{0}\in\mathcal{P}$ has no induced $C_{4}$ (Lemma
\ref{Lem_P_basic}(b)). Then Lemma \ref{Lem_C1_RL} and heredity, applied
repeatedly, give a connected member of $\mathcal{P}_{15}$ with deficiency at
least $3.$ This contradicts (Cens15). By Lemma \ref{Lem_M_Ftype}, $Y=aK_{1}\cup bK_{2}\cup tK_{3}$ with $t\leq1.$
The clique $S$ meets only one component of $Y.$ The isolated vertices of $Y$
outside $S$ are false twins with neighbourhood $\{o\}.$ So there is at most
one of them, $a\leq2$ and $b\geq(q-7)/2.$ Consider the components $K_{2}$ of
$Y$ disjoint from $S.$ There are at least $b-1\geq(q-9)/2$ of them, and each
is a pair of adjacent vertices with neighbourhood $\{o\}.$ We delete them one
at a time until the main component has order $13$ or $14.$ This gives a graph
$G^{\prime}=G_{0}^{\prime}\cup Z$ with $T(G)\geq T(G^{\prime})$ by Lemma
\ref{Lem_P_mono}. Here $G_{0}^{\prime}$ is connected and admissible, and
$\lambda_{1}(G_{0}^{\prime})\geq2$ as $o$ has at least $11$ neighbours. Also,
$G_{0}^{\prime}$ has no dominating vertex. Indeed, such a vertex would be
adjacent to $o$ and $w,$ i.e.\ it would lie in $S,$ and it would have at most
two neighbours in $Y.$ So $G^{\prime}$ belongs to the list of the case
$q\leq14,$ and $T(G)\geq T(G^{\prime})\geq2.$
\end{proof}

\bigskip\noindent\textbf{Concluding remarks.}~~The largest computations in our
proof are the census behind (Cens15) and the lists of admissible graphs of
order at most $14$ (Appendix \ref{App_comp}). A computer-free proof of
(Cens14) and (Cens15) would be of interest.

\appendix

\section{Computations}\label{App_comp}

\paragraph{Lists of admissible graphs, (Cens14) and (Cens15).}\label{App_lists}The
lists of members of $\mathcal{P}_{q},$ $q\leq15,$ were generated as described
in Section \ref{Sec_M}. We used independently written programs, and their
results agree. For $q\leq10$ the lists coincide with those obtained by brute
force over all graphs. The lists are \emph{inclusive}: we keep all graphs
whose computed cost is below $1+10^{-7}.$ So floating point errors can only
add graphs to them. For $q=2,\ldots,14$ they have
$2,4,11,33,104,267,545,857,1167,1454,1743,2032,2321$ members. (Cens14) and
(Cens15) hold for all graphs in these lists, so no graph of cost exactly $1$
has to be removed. The $49$ connected members of order $14$ with deficiency
$2$ are F-type with respect to every vertex of degree $12.$ No connected
member of order $15$ without induced $C_{4}$ (Lemma \ref{Lem_P_basic}(b)) has
deficiency at least $3.$

(Cens15) was also verified independently of these lists. By Lemma
\ref{Lem_P_basic}(b), Lemma \ref{Lem_C1_RL} and heredity, for $q\geq10$ every
connected member of $\mathcal{P}_{q}$ with deficiency at least $3$ arises from
such a member of order $q-1$ by adding one vertex. We started from a
brute-force computation for $q=9$ and extended step by step; there is no such
member of order $15.$ The details, including exact certificates for the graphs
of cost exactly $1$ met on the way, are given in the repository.

\paragraph{Small main components.}\label{App_small}We take every connected
graph $G_{0}$ of order $3\leq q\leq14$ in these lists with
$\operatorname{def}(G_{0})\geq2$ and computed $\lambda_{1}(G_{0})\geq
2-10^{-12}.$ For each such $G_{0}$ and every $Z$ of Lemma \ref{Lem_C1_struct}
with $q+\left\vert Z\right\vert $ even, we computed $T(G_{0}\cup Z)$ ($2944$
pairs; $4,56,120,496,296,544,260,400,192,376,200$ for $q=4,\ldots,14$). The
minimum, $2,$ is attained only
for $G_{0}=C_{6},$ $Z=\emptyset.$ It is present because the lists are
inclusive: $\overline{C_{6}}$ is the prism, with spectrum
$\{3,1,0^{2},(-2)^{2}\},$ so $\operatorname{cost}(C_{6})=1.$ The matrix
$J-A(C_{6})$ has the eigenvalues $4,2,1^{2},(-1)^{2},$ so $T(C_{6})=2.$ All
other values are at least $2.19.$

\paragraph{Exact computations and margins.}The following computations are
exact: the identities in rational arithmetic in the proof of Lemma
\ref{Lem_C1_order} and the value $T(C_{6})=2$ for small main components,
which follows from the spectrum given there. All other
computations are in floating point. They occur in the proof of Theorem
\ref{Tm_main} (two values of $T$), in Lemma \ref{Lem_C1_struct} (two costs and
three values of $T$), in the proof of Theorem \ref{Tm_cone} (four values of
$T$ and ten matrices $D^{\prime}(Z)$), in (Cens14), (Cens15) and for small
main components. They concern eigenvalues of symmetric matrices $A$ of order
$N\leq30$ with entries in $\{-1,0,1\}.$ By the backward stability of the
symmetric eigenvalue algorithms and Weyl's inequality, the computed eigenvalues
are correct up to a small multiple of $N\Vert A\Vert_{2}2^{-53}\leq
N^{2}2^{-53}<10^{-12}.$ Apart from the exact cases above, all margins are at
least $10^{-2},$ except for costs close to $1$ in the lists of
$\mathcal{P}_{q},$ $q\leq15;$ for this reason these lists are inclusive. The programs and the lists of graphs used here
will be provided with a link in the next version of the manuscript.

{\sloppy\bigskip\noindent\textbf{Acknowledgments.}~~The author \v{S}krekovski R.
acknowledges the partial support by ARIS program P1-0383, bilateral
Slovenian-Croatian project BI-HR/25-27-004 and the annual work program of
Rudolfovo, while the author Sedlar J. acknowledges the support by Project
KK.\allowbreak01.\allowbreak 1.\allowbreak1.\allowbreak02.\allowbreak0027
co-financed by the European Regional Development Fund, by Croatian Ministry of
Science, Education and Youth through the bilateral Croatian-Slovenian project
2025-26, and by the NextGeneration EU foundation via IP-UNIST-17 (GEORAZ). The
computations reported here were carried out in SageMath and in nauty. }

\bigskip

\noindent\textbf{Declaration of generative AI use.}~~During the preparation of this work the authors used Claude Opus 5.5 (Anthropic) to assist with editing and preparing the final version of the manuscript, including verification of calculations. The authors reviewed and edited all content and take full responsibility for the paper.

\end{document}